\catcode`\@=11

\font\tenmsa=msam10
\font\sevenmsa=msam7
\font\fivemsa=msam5
\font\tenmsb=msbm10
\font\sevenmsb=msbm7
\font\fivemsb=msbm5
\newfam\msafam
\newfam\msbfam
\textfont\msafam=\tenmsa  \scriptfont\msafam=\sevenmsa
  \scriptscriptfont\msafam=\fivemsa
\textfont\msbfam=\tenmsb  \scriptfont\msbfam=\sevenmsb
  \scriptscriptfont\msbfam=\fivemsb

\def\hexnumber@#1{\ifnum#1<10 \number#1\else
 \ifnum#1=10 A\else\ifnum#1=11 B\else\ifnum#1=12 C\else
 \ifnum#1=13 D\else\ifnum#1=14 E\else\ifnum#1=15 F\fi\fi\fi\fi\fi\fi\fi}

\def\msa@{\hexnumber@\msafam}
\def\msb@{\hexnumber@\msbfam}
\mathchardef\boxdot="2\msa@00
\mathchardef\boxplus="2\msa@01
\mathchardef\boxtimes="2\msa@02
\mathchardef\square="0\msa@03
\mathchardef\blacksquare="0\msa@04
\mathchardef\centerdot="2\msa@05
\mathchardef\lozenge="0\msa@06
\mathchardef\blacklozenge="0\msa@07
\mathchardef\circlearrowright="3\msa@08
\mathchardef\circlearrowleft="3\msa@09
\mathchardef\rightleftharpoons="3\msa@0A
\mathchardef\leftrightharpoons="3\msa@0B
\mathchardef\boxminus="2\msa@0C
\mathchardef\Vdash="3\msa@0D
\mathchardef\Vvdash="3\msa@0E
\mathchardef\vDash="3\msa@0F
\mathchardef\twoheadrightarrow="3\msa@10
\mathchardef\twoheadleftarrow="3\msa@11
\mathchardef\leftleftarrows="3\msa@12
\mathchardef\rightrightarrows="3\msa@13
\mathchardef\upuparrows="3\msa@14
\mathchardef\downdownarrows="3\msa@15
\mathchardef\upharpoonright="3\msa@16

\mathchardef\downharpoonright="3\msa@17
\mathchardef\upharpoonleft="3\msa@18
\mathchardef\downharpoonleft="3\msa@19
\mathchardef\rightarrowtail="3\msa@1A
\mathchardef\leftarrowtail="3\msa@1B
\mathchardef\leftrightarrows="3\msa@1C
\mathchardef\rightleftarrows="3\msa@1D
\mathchardef\Lsh="3\msa@1E
\mathchardef\Rsh="3\msa@1F
\mathchardef\rightsquigarrow="3\msa@20
\mathchardef\leftrightsquigarrow="3\msa@21
\mathchardef\looparrowleft="3\msa@22
\mathchardef\looparrowright="3\msa@23
\mathchardef\circeq="3\msa@24
\mathchardef\succsim="3\msa@25
\mathchardef\gtrsim="3\msa@26
\mathchardef\gtrapprox="3\msa@27
\mathchardef\multimap="3\msa@28
\mathchardef\therefore="3\msa@29
\mathchardef\because="3\msa@2A
\mathchardef\doteqdot="3\msa@2B

\mathchardef\triangleq="3\msa@2C
\mathchardef\precsim="3\msa@2D
\mathchardef\lesssim="3\msa@2E
\mathchardef\lessapprox="3\msa@2F
\mathchardef\eqslantless="3\msa@30
\mathchardef\eqslantgtr="3\msa@31
\mathchardef\curlyeqprec="3\msa@32
\mathchardef\curlyeqsucc="3\msa@33
\mathchardef\preccurlyeq="3\msa@34
\mathchardef\leqq="3\msa@35
\mathchardef\leqslant="3\msa@36
\mathchardef\lessgtr="3\msa@37
\mathchardef\backprime="0\msa@38
\mathchardef\risingdotseq="3\msa@3A
\mathchardef\fallingdotseq="3\msa@3B
\mathchardef\succcurlyeq="3\msa@3C
\mathchardef\geqq="3\msa@3D
\mathchardef\geqslant="3\msa@3E
\mathchardef\gtrless="3\msa@3F
\mathchardef\sqsubset="3\msa@40
\mathchardef\sqsupset="3\msa@41
\mathchardef\trianglerighteq="3\msa@44
\mathchardef\trianglelefteq="3\msa@45
\mathchardef\bigstar="0\msa@46
\mathchardef\between="3\msa@47
\mathchardef\blacktriangledown="0\msa@48
\mathchardef\blacktriangleright="3\msa@49
\mathchardef\blacktriangleleft="3\msa@4A
\mathchardef\blacktriangle="0\msa@4E
\mathchardef\triangledown="0\msa@4F
\mathchardef\eqcirc="3\msa@50
\mathchardef\lesseqgtr="3\msa@51
\mathchardef\gtreqless="3\msa@52
\mathchardef\lesseqqgtr="3\msa@53
\mathchardef\gtreqqless="3\msa@54
\mathchardef\Rrightarrow="3\msa@56
\mathchardef\Lleftarrow="3\msa@57
\mathchardef\veebar="2\msa@59
\mathchardef\barwedge="2\msa@5A
\mathchardef\doublebarwedge="2\msa@5B
\mathchardef\angle="0\msa@5C
\mathchardef\measuredangle="0\msa@5D
\mathchardef\sphericalangle="0\msa@5E
\mathchardef\varpropto="3\msa@5F
\mathchardef\smallsmile="3\msa@60
\mathchardef\smallfrown="3\msa@61
\mathchardef\Subset="3\msa@62
\mathchardef\Supset="3\msa@63
\mathchardef\Cup="2\msa@64

\mathchardef\Cap="2\msa@65

\mathchardef\curlywedge="2\msa@66
\mathchardef\curlyvee="2\msa@67
\mathchardef\leftthreetimes="2\msa@68
\mathchardef\rightthreetimes="2\msa@69
\mathchardef\subseteqq="3\msa@6A
\mathchardef\supseteqq="3\msa@6B
\mathchardef\bumpeq="3\msa@6C
\mathchardef\Bumpeq="3\msa@6D
\mathchardef\lll="3\msa@6E

\mathchardef\ggg="3\msa@6F

\mathchardef\circledS="0\msa@73
\mathchardef\pitchfork="3\msa@74
\mathchardef\dotplus="2\msa@75
\mathchardef\backsim="3\msa@76
\mathchardef\backsimeq="3\msa@77
\mathchardef\complement="0\msa@7B
\mathchardef\intercal="2\msa@7C
\mathchardef\circledcirc="2\msa@7D
\mathchardef\circledast="2\msa@7E
\mathchardef\circleddash="2\msa@7F
\def\ulcorner{\delimiter"4\msa@70\msa@70 }
\def\urcorner{\delimiter"5\msa@71\msa@71 }
\def\llcorner{\delimiter"4\msa@78\msa@78 }
\def\lrcorner{\delimiter"5\msa@79\msa@79 }
\def\yen{\mathhexbox\msa@55 }
\def\checkmark{\mathhexbox\msa@58 }
\def\circledR{\mathhexbox\msa@72 }
\def\maltese{\mathhexbox\msa@7A }
\mathchardef\lvertneqq="3\msb@00
\mathchardef\gvertneqq="3\msb@01
\mathchardef\nleq="3\msb@02
\mathchardef\ngeq="3\msb@03
\mathchardef\nless="3\msb@04
\mathchardef\ngtr="3\msb@05
\mathchardef\nprec="3\msb@06
\mathchardef\nsucc="3\msb@07
\mathchardef\lneqq="3\msb@08
\mathchardef\gneqq="3\msb@09
\mathchardef\nleqslant="3\msb@0A
\mathchardef\ngeqslant="3\msb@0B
\mathchardef\lneq="3\msb@0C
\mathchardef\gneq="3\msb@0D
\mathchardef\npreceq="3\msb@0E
\mathchardef\nsucceq="3\msb@0F
\mathchardef\precnsim="3\msb@10
\mathchardef\succnsim="3\msb@11
\mathchardef\lnsim="3\msb@12
\mathchardef\gnsim="3\msb@13
\mathchardef\nleqq="3\msb@14
\mathchardef\ngeqq="3\msb@15
\mathchardef\precneqq="3\msb@16
\mathchardef\succneqq="3\msb@17
\mathchardef\precnapprox="3\msb@18
\mathchardef\succnapprox="3\msb@19
\mathchardef\lnapprox="3\msb@1A
\mathchardef\gnapprox="3\msb@1B
\mathchardef\nsim="3\msb@1C
\mathchardef\napprox="3\msb@1D
\mathchardef\nsubseteqq="3\msb@22
\mathchardef\nsupseteqq="3\msb@23
\mathchardef\subsetneqq="3\msb@24
\mathchardef\supsetneqq="3\msb@25
\mathchardef\subsetneq="3\msb@28
\mathchardef\supsetneq="3\msb@29
\mathchardef\nsubseteq="3\msb@2A
\mathchardef\nsupseteq="3\msb@2B
\mathchardef\nparallel="3\msb@2C
\mathchardef\nmid="3\msb@2D
\mathchardef\nshortmid="3\msb@2E
\mathchardef\nshortparallel="3\msb@2F
\mathchardef\nvdash="3\msb@30
\mathchardef\nVdash="3\msb@31
\mathchardef\nvDash="3\msb@32
\mathchardef\nVDash="3\msb@33
\mathchardef\ntrianglerighteq="3\msb@34
\mathchardef\ntrianglelefteq="3\msb@35
\mathchardef\ntriangleleft="3\msb@36
\mathchardef\ntriangleright="3\msb@37
\mathchardef\nleftarrow="3\msb@38
\mathchardef\nrightarrow="3\msb@39
\mathchardef\nLeftarrow="3\msb@3A
\mathchardef\nRightarrow="3\msb@3B
\mathchardef\nLeftrightarrow="3\msb@3C
\mathchardef\nleftrightarrow="3\msb@3D
\mathchardef\divideontimes="2\msb@3E
\mathchardef\varnothing="0\msb@3F
\mathchardef\nexists="0\msb@40
\mathchardef\mho="0\msb@66
\mathchardef\thorn="0\msb@67
\mathchardef\beth="0\msb@69
\mathchardef\gimel="0\msb@6A
\mathchardef\daleth="0\msb@6B
\mathchardef\lessdot="3\msb@6C
\mathchardef\gtrdot="3\msb@6D
\mathchardef\ltimes="2\msb@6E
\mathchardef\rtimes="2\msb@6F
\mathchardef\shortmid="3\msb@70
\mathchardef\shortparallel="3\msb@71
\mathchardef\smallsetminus="2\msb@72
\mathchardef\thicksim="3\msb@73
\mathchardef\thickapprox="3\msb@74
\mathchardef\approxeq="3\msb@75
\mathchardef\succapprox="3\msb@76
\mathchardef\precapprox="3\msb@77
\mathchardef\curvearrowleft="3\msb@78
\mathchardef\curvearrowright="3\msb@79
\mathchardef\digamma="0\msb@7A
\mathchardef\varkappa="0\msb@7B
\mathchardef\hslash="0\msb@7D
\mathchardef\hbar="0\msb@7E
\mathchardef\backepsilon="3\msb@7F
\def\Bbb{\ifmmode\let\next\Bbb@\else
 \def\next{\errmessage{Use \string\Bbb\space only in math mode}}\fi\next}
\def\Bbb@#1{{\Bbb@@{#1}}}
\def\Bbb@@#1{\fam\msbfam#1}

\catcode`\@=\active

\def\sw#1{{\sb{(#1)}}}
\def\bC{{\overline{C}}}
\def\sq{{S\sp 2\sb q(\mu ,\nu)}}
\def\rco#1{{\Delta\sb{#1}}}
\def\lco#1{{~\sb{#1}\!\Delta}}

\documentstyle[12pt]{article}
\normalbaselines
\catcode `\@=11
\@addtoreset{equation}{section}

\def\section{\@startsection {section}{1}{\z@}{-3.5ex plus -1ex minus
     -.2ex}{2.3ex plus .2ex}{\normalsize\bf}}
\def\subsection{\@startsection{subsection}{2}{\z@}{-3.25ex plus -1ex minus
 -.2ex}{1.5ex plus .2ex}{\normalsize\bf}}

\def\thebibliography#1{\section*{References\markboth
  {REFERENCES}{REFERENCES}}\list
  {[\arabic{enumi}]}{\settowidth\labelwidth{[#1]}\leftmargin\labelwidth
  \advance\leftmargin\labelsep
  \usecounter{enumi}}
  \def\newblock{\hskip .11em plus .33em minus -.07em}
  \sloppy
  \sfcode`\.=1000\relax}
 
\catcode `\@=12

\begin{document}
\baselineskip 22pt
{\ }\qquad\qquad \hskip 4.0in DAMTP/97-18

\vspace*{2.5cm}
\noindent
{ \bf QUANTUM HOMOGENEOUS SPACES AND
  COALGEBRA BUNDLES}\vspace{1.3cm}\\
\noindent
\hspace*{1in}
\begin{minipage}{13cm}
Tomasz Brzezi\'nski\vspace{0.3cm}\\ 
DAMTP, University of Cambridge, Cambridge CB3 9EW, U.K.\\
\end{minipage}

\vspace*{0.5cm}

\begin{abstract}
\noindent
It is shown that quantum homogeneous spaces of a quantum group $H$ can
be viewed as fibres of quantum fibrations  with the total space $H$
that are dual to  coalgebra bundles. As concrete examples
of such structures the 
fibrations with the
quantum 2-sphere and the quantum hyperboloid fibres are considered.
\end{abstract}

\section{\hspace{-4mm}.\hspace{2mm}INTRODUCTION}
Among the variety of quantum spaces quantum homogeneous spaces seem to
play a special role. Their symmetry structure is rich enough to display
various properties of quantum group actions. Recently by studying the
structure of quantum homogeneous spaces 
\cite{Brz:hom} we were lead to the notion of a coalgebra principal
bundle \cite{BrzMa:coa} which generalises the notion of a quantum
principal bundle \cite{BrzMa:gau} \cite{Dur:geo} \cite{Pfl:fib} or a
Hopf-Galois extension 
\cite{KreTak:hop} \cite{Mon:hop}. A coalgebra principal bundle is
built with a coalgebra $C$ which plays the role of the `structure
group' (fibre) and two algebras $M$ and $P$, the first of which plays
the role of the base manifold and the second corresponds to the total
space. The spaces $C$, $M$, $P$ satisfy some conditions which in
algebraic terms  mean that $P$ is an extension of $M$
by $C$, known as a C-Galois extension (see Section~2 for 
definition). 
Since all the spaces involved in the definition of coalgebra principal
bundles are either algebras or coalgebras and  the
notions of a coalgebra and an algebra are dual to each other, it seems
natural to consider an object dual to a coalgebra principal
bundle. This is built with an algebra $A$ which plays the role of a fibre
and the coalgebras $\bC$ and $C$, the first of which is a base and the
second is a total space. Again, in algebraic terms $C$ is an
extension of the coalgebra $\bC$ by an algebra $A$ which can be termed
an {\em A-Galois coextension} (see Section~3 for definition). In this
paper we show that the structure 
of quantum
homogeneous spaces provides natural  examples of such A-Galois coextensions. 

Quantum homogeneous spaces $M$ of Hopf algebras $H$ that we
study in Section~4 have the following remarkable property. 
As an algebra $H$  is  a C-Galois extension of $M$ by a
certain coalgebra $C$, while  as a coalgebra it is an A-Galois
coextension 
of (the same) $C$ by $M$\footnote{We would like to stress, however,
that from the point of view of 
entwining structures involved in  the definition of coalgebra bundles in
\cite{BrzMa:coa} those two bundles differ substantially. In the first case the
spaces $C$ and $H$  while in the second $H$ and $M$ are entwined.}.  
 In geometric terms  these quantum homogeneous
spaces play the role of a base manifold in the former case and a fibre
in the latter (the same total space in both cases). This suggests  the
notion of a CA-Galois biextension, 
which is very much reminiscent of
bicrossproducts of \cite{Ma:phy}, and we think might be well-worth
studying. 

Throughout the paper, all vector spaces are over a field $k$ of
characteristic 0, although the results can be extended to more
general fields or even commutative rings. By an algebra we mean
an associative algebra over $k$ with unit 
denoted by 1. In a
coalgebra $C$, $\Delta$ 
is the coproduct and 
$\epsilon :C\to k$ is the counit. We use the Sweedler notation for a
coproduct, $\Delta(c) = c\sw 1\otimes c\sw 2$ (summation 
understood), for any $c\in C$.

\section{\hspace{-4mm}.\hspace{2mm}COALGEBRA BUNDLES AND EMBEDDABLE H-SPACES }
A coalgebra principal bundle is
defined as follows. Let $C$ be a coalgebra, $P$ be an algebra and a right
$C$-comodule. It means that there is a linear map  $\Delta_P : P\to P\otimes
C$, such that $(id\otimes \Delta)\circ\rco P = (\rco P\otimes
id)\circ\rco P$ and $(id\otimes\epsilon)\circ\rco P = id$. The map
$\rco P$ is known as a {\em right coaction}. Let $e\in C$ be
group-like, i.e. $\Delta e = e\otimes e$, and let $M=P_e^{coC} = \{x\in
P | \Delta_P (x) = x\otimes e\}$ be a space of coinvariants. Furthermore
assume that the coaction is 
left-linear over $M$, i.e. $\Delta_P(xu) = x\Delta_P (u)$ for all $x\in
M$, $u\in P$ and that $\Delta_P (1) = 1\otimes e$. Then $M$ is an
algebra and the canonical
map
\begin{equation}
\chi_M:  P\otimes_M P \to P\otimes C, \qquad \chi_M: u\otimes v
\mapsto u\Delta_P(v)
\label{chim}
\end{equation}
is well defined. We say that $P$ is a {\em $C$-coalgebra principal bundle
over $M$} or a {\em C-Galois extension of $M$} if the map $\chi_M$ is
a bijection\footnote{In 
\cite{BrzMa:coa} a coalgebra principal bundle was defined in the
framework of entwining structures. It is shown in \cite{BrzHaj:not}
that the above definition is equivalent to the one in
\cite{BrzMa:coa}. The same remark applies to the definition of a dual
coalgebra bundle in Section~3.}. We denote this bundle by 
$P(M,C,e)$. This definition reduces to  the definition of a quantum
principal bundle with universal differential calculus of
\cite{BrzMa:gau} if $C$ is a Hopf algebra, 
$\Delta_P$ is an algebra map and $e=1\in C$
(cf. \cite[Lemma~3.2]{Brz:tra} \cite[Proposition~1.6]{Haj:str}). 

The notion of a coalgebra
bundle is motivated by the structure of left quantum homogeneous
spaces and, indeed, this is the context in which $C$-Galois extensions
appeared for the first time in \cite{Sch:nor} (although studied for
entirely different reasons there). Let $H$ be
a Hopf algebra and $M$ be a left (resp. right) $H$-comodule
algebra. It means that there is a linear map  $\lco M : M\to
H\otimes M$ such that $(id\otimes \lco M)\circ\lco M = (\Delta\otimes
id)\circ\lco M$ and $(\epsilon\otimes id)\circ\lco M = id$, i.e. a
left coaction (resp. a right coaction $\Delta_M : M\to M\otimes H$)
which, furthermore, is an algebra map. Following \cite{Pod:sym} we say
that $M$ is a
{\it left ({\rm resp.} right) embeddable quantum
homogeneous space} or simply a {\it left ({\rm resp.} right) embeddable
H-space} if there exists an algebra inclusion
$i: M\hookrightarrow H$ such that $\Delta\circ i
= (id\otimes i)\circ\lco M$ ( resp. $\Delta\circ i
= (i\otimes id)\circ\Delta\sb M$). When $M$ is viewed as a subalgebra
of $H$ via $i$ then 
the coaction coincides with the coproduct in $H$.

Given a left embeddable $H$-space $M$ one  defines the right
ideal $J_R = i(M)^+ H$,
where $i(M)^+ = \ker\epsilon\cap i(M)$.
This ideal is a coideal of $H$,
i.e. $\Delta(J_R)\subset H\otimes J_R \oplus J_R\otimes H$ and,
obviously, $\epsilon(J_R) =0$. Therefore $C_R = H/J_R$ is a coalgebra
and a canonical surjection $\pi_R : H \to C_R = H/J_R$ is a coalgebra
map. Furthermore, $C_R$ is a right $H$-module with the action $\rho_R(c,h)
=\pi_R(gh)$ for any $h\in H$, $c\in C_R$ and
$g\in\pi_R^{-1}(c)$. Finally $H$ is a 
right $C$-comodule with the coaction $\Delta_H = (id\otimes
\pi_R)\circ\Delta$. If one chooses $e=\pi_R(1)$ then $H_e^{co C_R}$ is a
subalgebra of $H$ and $\Delta_H$ is left-linear over $H_e^{co
C_R}$.  Also, the canonical map $\chi_{H_e^{co
C_R}}$ is a bijection \cite[Lemma~1.3]{Sch:nor} or
\cite[Example~4.4]{Brz:hom}. Therefore $H$ is a $C$-principal
coalgebra bundle over $H_e^{co
C_R}$. Furthermore,  $i(M) \subset H_e^{co C_R}$, i.e.  $M$ is a subalgebra of
$H_e^{co C_R}$. If $M$ is 
isomorphic to $H_e^{co C_R}$ via $i$ then clearly $H$ is a coalgebra
bundle over $M$. For example, the map $i : M \hookrightarrow H_e^{co
C_R}$ is an 
isomorphism, if $H$ is 
faithfully flat right or left $M$-module \cite[Lemma~1.3]{Sch:nor} (cf.
\cite[Section~I.3]{Bou:com} for the definition and discussion of
faithful flatness). 

\section{\hspace{-4mm}.\hspace{2mm}DUAL
COALGEBRA BUNDLES AND EMBEDDABLE H-SPACES} 
First recall the definition of a dual coalgebra bundle 
\cite{BrzMa:coa} which generalises the object analysed in
\cite[Theorem~II]{Sch:pri}.  
Let $A$ be an algebra, let $C$ be a coalgebra and a right $A$-module
with the action $\rho_C : C\otimes A\to C$ and let $\kappa: A\to k$ be
an algebra character. Define a vector space  $J_\kappa =
\rho_C(C,\ker\kappa)$ 
and a quotient space
$\bC =  C/J_\kappa$. Let $\pi_\kappa : C \to \bC$ be a
canonical surjection. It can be easily shown \cite{BrzHaj:not} that if
$\epsilon\circ\rho_C = \epsilon\otimes\kappa $ and 
\begin{equation}
\Delta_\kappa\circ\rho_C =
( id \otimes\rho_C)\circ(\Delta_\kappa\otimes id), 
\label{colinear}
\end{equation}
where $\Delta_\kappa = (\pi_\kappa\otimes id)\circ\Delta$, then
$J_\kappa$ is a coideal and therefore $\bC$ is a coalgebra. 
Since $\bC$ is a coalgebra and $\pi_\kappa$ is a
coalgebra map, $C$ is a left $\bC$ comodule with the coaction
$\Delta_\kappa$. Since, moreover, $C$ is a right $A$-module, 
(\ref{colinear}) is 
equivalent to the right-linearity of the coaction $\Delta_\kappa$ over
$A$. Next, one defines the cotensor product by
$$
C\square_\bC C = \{b\otimes c\in C\otimes C | b\sw 1\otimes
\pi_\kappa(b\sw 2)\otimes c = 
b\otimes \pi_\kappa(c\sw 1)\otimes c\sw 2\}.
$$
Define a canonical map
\begin{equation}
\chi^\bC :C\otimes A\to C\square\sb\bC C , \qquad   \chi ^\bC:  c\otimes
 a \mapsto c\sw 
 1\otimes \rho_C(c\sw 2, a).
\label{chim2}
\end{equation}
We say
that $C$ is a {\em dual coalgebra $A$-principal bundle over $\bC$} or
an {\em A-Galois coextension of $\bC$} if
the canonical 
map $\chi^\bC$ is a bijection. This bundle is denoted by $C(\bC, A,\kappa)$.

We show that quantum embeddable homogeneous spaces provide
examples of dual coalgebra bundles. Consider a {\em right} $H$-embeddable
space $M$ with embedding $i:M\hookrightarrow H$. Similarly as for left
spaces, $M$ can be interpreted as a quantum
quotient space, i.e. as coinvariants of the coaction of certain 
coalgebra. Consider a left 
ideal $J_L = Hi(M)^+$.
This ideal is a coideal of $H$,
therefore $C_L = H/J_L$ is a coalgebra
and a canonical surjection $\pi_L : H \to C_L$ is a coalgebra
map. $C_L$ is a left $H$-module with the action $\rho_L(h,c)
=\pi_L(hg)$ for any $h\in H$, $c\in C_L$ and
$g\in\pi_L^{-1}(c)$. Finally $H$ is a 
left $C$-comodule with the coaction $\lco H = (\pi_L\otimes
id)\circ\Delta$. If one chooses $e=\pi_L(1)$ then the space of
coinvariants $^{co C_L}H_e =
\{ x\in H | \lco H (x) = e\otimes x\}$ is a
subalgebra of $H$  and, in nice cases  such as the faithfully flat one,  $M$ is
isomorphic to $^{co 
C_L}H_e$ by $i$. From now on we assume that $i(M)=~^{co 
C_L}H_e$. $M$ acts on $H$ from the right via the
multiplication, i.e. $\rho_H(h,x) = hi(x)$, 
thus $J_L$ has the structure of $J_\kappa$ above with $\kappa  =
\epsilon\circ i$. Clearly $\epsilon\circ\rho_H =
\epsilon\otimes\kappa$. Moreover, $\lco H$ is right-linear over $M$,
hence (\ref{colinear}) holds. 
Thus we are in the setting
needed for a dual coalgebra bundle over
$C_L$, with a fibre $M$ and the total space $H$. It remains to be checked
that $\chi^{C_L}$ is a bijection. Consider a map 
$$
\tilde{\chi} : H\square\sb{C_L} H \to H\otimes M, \qquad
\tilde{\chi} : h\otimes g \mapsto h\sw 1\otimes i\sp{-1}((Sh\sw 2)g),
$$
where $S:H\to H$ is the antipode of $H$.
First it needs to be shown that  
$\tilde{\chi}$ is well-defined.
Take any $h\otimes g \in H\square\sb{C_L}
H$ and apply $id\otimes\lco H$ to $h\sw 1\otimes (Sh\sw 2)g$,
\begin{eqnarray*}
h\sw 1 \otimes \pi_L((Sh\sw
3)g\sw 1)\otimes (Sh\sw
2)g\sw 2 & = & h\sw 1\otimes \rho_L(Sh\sw 3, \pi_L(g\sw 1)) \otimes (Sh\sw
2)g\sw 2 \\
& = &h\sw 1\otimes \rho_L(Sh\sw 3, \pi_L(h\sw 4)) \otimes (Sh\sw
2)g \\
& = & h\sw 1\otimes e \otimes (Sh\sw
2)g.
\end{eqnarray*}
The second equality follows from  the fact that $h\otimes g \in
H\square\sb{C_L}H$. Therefore $h\sw 1\otimes (Sh\sw 2)g\in H\otimes
~^{co C_L}H_e$ and $\tilde{\chi}$ is well-defined. An elementary calculation
shows that $\tilde{\chi}$ is an inverse of $\chi^{C_L}$. Consequently
there is a dual coalgebra bundle $H(C_L, M, \epsilon\circ i)$.

It is perhaps worth mentioning  that every left homogeneous $H$-space
$M$ is a right $H^{cop}$-space, where $H^{cop}$ is isomorphic to $H$
as an algebra but has an opposite coproduct.
If $H$ has the invertible antipode then $H^{cop}$ is also a Hopf
algebra. Viewing a  left embeddable $H$-space $M$ as a right
$H^{cop}$-space one 
can thus construct the bundle $H^{cop}(C_L^{cop},M,\epsilon\circ
i)$. Therefore $M$ plays a double role: firstly it is a base for
$H(M,C_R,\pi_R(1))$ constructed in Section~2 and secondly it is a fibre
of $H^{cop}(C_L^{cop},M,\epsilon\circ
i)$. 
\section{\hspace{-4mm}.\hspace{2mm}TWO EXPLICIT EXAMPLES}
The examples below  are selected from a long list of
well-known quantum homogeneous spaces such as quantum planes,
homogeneous spaces of $E_\kappa(2)$ and $E_\kappa(3)$ etc, to all of which
the above construction applies too. The first example is a
non-trivial quantum spherical fibration of $SU_q(2)$ while the second one
is a trivial quantum hyperbolic fibration of $E_q(2)$ (the word
`trivial' here means that the total space is a certain crossed product
of a base and a fibre see \cite{BrzMa:coa} and \cite{Brz:cro} for
details and examples). 
\subsection{\hspace{-5mm}.\hspace{2mm}Quantum Spherical Fibration of
$SU_q(2)$} 
The quantum two-sphere \cite{Pod:sph} $S\sb q \sp 2(\mu ,\nu)$ is
a polynomial complex algebra generated by the unit and $x$, $y$, $z$,
and the relations
$$
 xz = q\sp 2zx,  \qquad xy  = -q(\mu -z)(\nu +z), \quad 
 yz  = q\sp{-2}zy, \qquad yx  = -q(\mu -q\sp{-2}z)(\nu +q\sp{-2}
 z),
$$
where $\mu, \nu $ and $q\neq 0$ are real parameters, $\mu\nu \geq 0$,
$(\mu ,\nu) \neq (0,0)$.
The quantum sphere is a $*$-algebra with the
$*$-structure $x\sp * = -q y$, $z\sp * =z$.

The quantum sphere $S\sp 2\sb q(\mu ,\nu)$ is an $SU\sb q(2)$
homogeneous quantum space. $SU\sb q(2)$ is a complex algebra generated
by the matrix of 
generators
${\bf t} = \pmatrix{\alpha &\beta \cr \gamma &\delta}$
and the
relations
$$
\alpha\beta =q\beta\alpha, \quad \alpha\gamma =q\gamma\alpha, \quad
\alpha\delta = \delta\alpha + (q-q\sp{-1})\beta\gamma ,
$$
$$\beta\gamma = \gamma\beta, \quad \beta\delta = \delta\beta ,
\quad \gamma\delta=q\delta\gamma, \quad \alpha\delta - q\beta\gamma =1 .
$$
$SU\sb q(2)$ is a Hopf algebra of a matrix group type, i.e.
$$
\Delta {\bf t} = {\bf t}\otimes{\bf t}, \quad
\epsilon{\bf t} = 1,
\quad S{\bf t} = \pmatrix{ \delta &-q\sp{-1}\beta \cr -q\gamma &\alpha},
$$
and has a $*$-structure given
by $\delta = \alpha\sp *$, $\gamma = -q\sp{-1} \beta\sp *$. 
$\sq$ is not only a quantum homogeneous
space but also it is an embeddable $SU\sb q(2)$-space. Furthermore it
is a left and right embeddable $SU\sb q(2)$-space \cite{Pod:sym} with
the $*$-algebra  inclusions
$$
i\sb L(x) = \sqrt{\mu\nu}(q\alpha\sp 2 - \beta\sp 2) +
(\mu -\nu)\alpha\beta, \quad
i\sb L(y) = \sqrt{\mu\nu}(q\gamma\sp 2 - \delta\sp 2) +
(\mu -\nu)\gamma\delta,
$$
$$
i\sb L(z) = -\sqrt{\mu\nu}(q\alpha\gamma - \beta\delta) -
q(\mu -\nu)\beta\gamma ,
$$
$$
i\sb R(x) = q\sqrt{\mu\nu}(\alpha\sp 2 - q\gamma\sp 2) -
q(\mu -\nu)\alpha\gamma, \quad
i\sb R(y) = \sqrt{\mu\nu}(q^{-1}\beta\sp 2 - \delta\sp 2) -
q^{-1}(\mu -\nu)\beta\delta,
$$
$$
i\sb R(z) = \sqrt{\mu\nu}(\alpha\beta - q\gamma\delta) -
q(\mu -\nu)\beta\gamma .
$$
To simplify the forthcoming analysis assume that
$\mu\neq\nu$. Then
$\sq$ depends on two real parameters $q$
and $p = \frac{\sqrt{\mu\nu}}{\mu-\nu}$. Applying the analysis of
Section~2 to $\sq$ and using the map $i_L$ above one defines a right
ideal $J_R$ in
$SU\sb q(2)$ generated by 
$$
p(q\alpha\sp 2 - \beta\sp 2) + \alpha\beta -pq, \quad
p(q\gamma\sp 2 - \delta\sp 2) + \gamma\delta + p, \quad
p(q\alpha\gamma - \beta\delta) + q\beta\gamma .
$$
and  constructs the coalgebra
$C_R = SU\sb q(2)/J_R$. $C_R$ is spanned by
$\pi_R(\alpha^n)$, $\pi_R(\delta^n)$, $n = 0,1,2,\ldots$ It is shown
in \cite{Brz:hom} that the space of coinvariants $M = {SU\sb q(2)}^{co
C_R}_e$ is isomorphic to 
$\sq$. Therefore  a coalgebra principal bundle
$SU_q(2)(\sq, C_R, \pi_R(1))$ over $\sq$ with fibre $C_R$ is constructed.

On the other hand we can apply the analysis of Section~3 and use the
map $i_R$ which allows us to view $\sq$ as a right embeddable
$SU\sb q(2)$-space. We proceed to define a left ideal $J_L$ in
$SU\sb q(2)$ generated by 
$$
p(q\delta\sp 2 - \beta\sp 2) + \beta\delta -pq, \quad
p(q\gamma\sp 2 - \alpha\sp 2) + \alpha\gamma + p, \quad
p(\alpha\beta - q\gamma\delta) - q\beta\gamma .
$$
and  construct the coalgebra
$C_L = SU\sb q(2)/J_L$. $C_L$ is spanned by 
$\pi_L(\alpha^n)$, $\pi_L(\delta^n)$, $n = 0,1,2,\ldots$.  Following
\cite{Brz:hom} one deduces that the 
 space of left coinvariants $M = ~^{co C_L}{SU\sb
q(2)}_e$ is isomorphic to 
$\sq$ and there is  a dual coalgebra bundle $SU_q(2)(C_L,\sq ,
\kappa)$ over $C_L$ with fibre $\sq$. Here $\kappa$ is a
$*$-algebra character of $\sq$ given by $\kappa(x) = q\sqrt{\mu\nu}$,
$\kappa(y) = -\sqrt{\mu\nu}$ and $\kappa(z) = 0$. 

Remarkably,
it can be shown that although  $C_L$ and $C_R$ are different as
$SU_q(2)$ modules they are isomorphic to each other as coalgebras. The
isomorphism is $\pi_L(\alpha^n) \mapsto \pi_R(\delta^n)$ and
$\pi_L(\delta^n) \mapsto \pi_R(\alpha^n)$. Thus we identify $C_L$
with $C_R$ as coalgebras and conclude that $SU\sb q(2)$ can be viewed as a
CA-Galois biextension: it is a C-Galois extension of $\sq$ by $C_R$
and it is an A-Galois coextension of $C_R$ by $\sq$. 

\subsection{\hspace{-5mm}.\hspace{2mm}Quantum Hyperbolic Fibration
of $E_q(2)$}
As another example take the quantum hyperboloid $X_q$
\cite{Sch:int}. $X_q$ is a complex algebra generated by $z_+$, $z_-$
and the identity, 
and the relation 
$$
z_+z_- = q^2 z_- z_+ +(1-q^2),
$$
where $q$ is a real number. $X_q$ can also be viewed as a deformation
of the complex plane obtained by the stereographic projection of
$S_q^2(\mu, 0)$ \cite{ChuHo:sph}. $X_q$ is a $*$-algebra with involution
$z_+^* = z_-$ and it
is  a quantum homogeneous space of  $E_q(2)$. The latter is a complex
algebra 
generated by $v,v^{-1}, n_+,n_-$ subject to the following  
relations \cite{VakKor:alg} \cite{Wor:unb}
$$
vn_\pm = q^2n_\pm v, \quad 
\quad n_+n_-
=q^2 n_-n_+ \quad vv^{-1} = v^{-1} v =1.
$$
$E_q(2)$ is a quantum group with a coproduct, counit and the antipode:
$$
\Delta v^{\pm 1} = v^{\pm 1}\otimes v^{\pm 1}, \quad \Delta n_\pm =
n_\pm \otimes 1 + v^{\pm 1}\otimes n_\pm,$$
$$
\epsilon(v^{\pm 1}) =
1, \quad \epsilon(n_\pm) = 0, \quad S(v^{\pm 1}) =v^{\mp 1}, \quad
S(n_\pm) = -v^{\mp 1}n_\pm .
$$
$E_q(2)$ is a $*$-algebra, $v^* =
v^{-1}$, $n_+^*=n_-$. It is shown in \cite{BonCic:fre} that $X_q$ is a
left embeddable 
$E_q(2)$-space with the coaction $~_{X_q}\Delta(z\sb\pm) =
n\sb\pm\otimes 1 +v\sp{\pm 1}\otimes z\sb\pm$ and the  $*$-algebra embedding
$i_L(z_\pm) = v^{\pm 1} + 
n_\pm$. It can be 
easily seen that $X_q$ is also a right 
embeddable $E_q(2)$-space with the coaction $\Delta_{X_q}(z\sb\pm) =
z\sb\pm\otimes v\sp{\pm 1} +q\sp{\mp 1}\otimes v\sp{\pm 1}n\sb\mp$ and the
$*$-algebra embedding $i_R(z_\pm) = v^{\pm 1} 
+q^{\mp 1}v^{\pm 1}n_\mp$. In both cases $X_q$ comes from the construction
described in Sections~2 and 3 with the right ideal $J_R$ generated by
$v^{\pm 1} + n_\pm -1$ and the left ideal $J_L$ generated by $v^{\pm 1}
-q\sp{\pm 1}n_\pm -1$. The corresponding coalgebras $C_R=E_q(2)/J_R$
and $C_L = E_q(2)/J_L$ are
isomorphic to each other as coalgebras (but, of course, differ by
their $E_q(2)$- module structure). $C_R$ is spanned by group-like $\pi_R(v^n)$,
$n\in \bf Z$ and $C_L$ is spanned by group-like $\pi_L(v^n)$,
$n\in \bf Z$, and the isomorphism is given by $\pi_R(v^n)\mapsto
\pi_L(v^n)$. Thus we can identify $C_L$, $C_R$ with the coalgebra $C$
spanned by elements $c_p$, $p\in \bf Z$, such that $\Delta c_p =
c_p\otimes c_p$. The
right and left actions of $E_q(2)$ on $C$ are given by 
$$
\rho_R(c_p, v^k) =  c_{p+k}, \quad \rho_R(c_p,n_\pm) =
q^{2p}(c_p -c_{p\pm 1}),
$$
$$
\rho_L(v^k,c_p) = c_{p+k}, \quad \rho_L(n_\pm,c_p) = q^{2p\mp 1}(c_{p\pm
1}-c_p).
$$
Therefore the quantum Euclidean group $E_q(2)$ is a
CA-Galois biextension built on $C$ and $X_q$. In other words we have
a coalgebra and a dual coalgebra principal bundles,
$E_q(2)(X_q,C,c_0)$ over $X_q$ with fibre $C$ and  $E_q(2)(C,X_q,
\kappa)$ over $C$ with fibre $X_q$. The map $\kappa$ is a
$*$-character of $X_q$ given by $\kappa(z\sb\pm) = 1$.

To obtain a slightly more topological insight into the structure of
$E_q(2)$ one equips $C$ with an algebra structure of ${\bf
C}[Z,Z^{-1}]$ by $c_n\mapsto Z^n$ and thus turns it into a Hopf
algebra. By requiring that $\pi_L$ and 
$\pi_R$ be $*$-maps the compatible $*$-structure is obtained on ${\bf
C}[Z,Z^{-1}]$ as $Z^* = Z^{-1}$. Thus $C\equiv {\bf
C}[Z,Z^{-1}]$ is an algebra of functions on a (classical)
circle $S^1$. $E_q(2)$ can be viewed then as the $S^1$-bundle over
$X_q$ as well as a dual bundle over a  circle with quantum hyperbolic
fibres.
\section*{\hspace{-5mm}~\hspace{2mm}ACKNOWLEDGEMENT}
Research supported by the EPSRC grant
GR/K02244.

\end{document}